\documentclass{ifacconf}

\usepackage{graphicx}      
\usepackage{natbib}        
\usepackage{amsmath,amsfonts,amssymb}

\begin{document}
\begin{frontmatter}

\title{Performance Bounds of Model Predictive Control for Unconstrained and Constrained Linear Quadratic Problems and Beyond\thanksref{footnoteinfo}} 

\thanks[footnoteinfo]{This work was supported by the Swedish Foundation for Strategic Research, the Swedish Research
Council, and the Knut and Alice Wallenberg Foundation. It is also supported by the Swiss National Science Foundation under NCCR Automation (grant agreement 51NF40\_180545).}

\author[First]{Yuchao Li}
\author[Second]{Aren Karapetyan} 
\author[Second]{John Lygeros} 
\author[First]{Karl H. Johansson} 
\author[First]{Jonas M\aa rtensson}

\address[First]{Division of Decision and Control Systems, KTH Royal Institute of Technology, Sweden, (e-mail: yuchao,kallej,jonas1@kth.se).}

\address[Second]{Automatic Control Laboratory, Swiss Federal Institute of Technology in Zürich, Switzerland, (e-mail: akarapetyan,jlygeros@ethz.ch).}

\begin{abstract}
We study unconstrained and constrained linear quadratic problems and investigate the suboptimality of the model predictive control (MPC) method applied to such problems. Considering MPC as an approximate scheme for solving the related fixed point equations, we derive performance bounds for the closed-loop system under MPC. Our analysis, as well as numerical examples, suggests new ways of choosing the terminal cost and terminal constraints, which are \emph{not} related to the solution of the Riccati equation of the original problem. The resulting method can have a larger feasible region, and cause hardly any loss of performance in terms of the closed-loop cost over an infinite horizon.
\end{abstract}


\end{frontmatter}

\section{Introduction}
Model predictive control (MPC) is a well-established scheme for constrained optimal control problems with continuous state and control spaces. In its most widely adopted form \citep{scokaert1998constrained} a system model is used to predict its performance over some finite number of stages, and a terminal cost is added to account for the trajectory beyond this prediction horizon. Moreover, an additional constraint is imposed on the state at which the prediction ends. Suitably designed terminal cost and constraint lead to desirable properties \citep{mayne2000constrained}. A recent monograph \citep{bertsekas2022lessons} shows that this structure of MPC is quite similar to the algorithms applied in the high-profile successes in the field of reinforcement learning \citep{silver2017mastering}. It has also developed a conceptual framework centered around dynamic programming (DP), which is useful for analyzing the performance of MPC. 

In this work, we apply the tools introduced in \citep{bertsekas2022lessons}, and investigate the impact of the terminal cost on the performance of MPC measured by the closed-loop cost accumulated over infinite stages. Using unconstrained and constrained linear quadratic regulation (LQR) problems as a vehicle, we derive performance bounds of MPC applied to those problems compared against optimal control. The insights gained from our analysis suggest new designs of those terminal ingredients, which likely make MPC feasible for a larger set of states, while cause little to no degradation in its performance.   

For the LQR problem, a performance bound of MPC is given \citep{bitmead1985monotonicity}. It relies on the monotonicity property, which holds well beyond the problem considered here, as discussed in \citep{bertsekas1975monotone,bertsekas1977monotone}. A suboptimality analysis that is similar in purpose to our study is reported in \citep{grune2008infinite}. It relies in part on a relaxed form of DP introduced in \citep{rantzer2006relaxed}. Our work makes different assumptions, and leads to new design choices of the terminal ingredients, and may be regarded as complementary to that of \citep{grune2008infinite}.

Unlike MPC where the performance bound analysis is relatively rare, bounds are pervasive in the study of Markovian decision problems (MDP), with some classical results presented in \citep{denardo1967contraction}. \citep{bertsekas2022lessons} takes an abstract approach that unifies the MDP and the problems addressed by MPC. It regards MPC as \emph{approximation in value space}, where the \emph{approximation} refers to the fact that terminal ingredients act together to approximate the real optimal cost. It also points out that the typical bounds in MDP may be too conservative to reflect the true performance, which we confirm by deriving better bounds for MPC.

Another concept that is related to our study is the investigation of \emph{regret} in dynamical systems. Regret is a non-asymptotic performance metric for optimal control problems involving uncertainty. Given a controller, regret is measured by comparing its associated cost over a finite number of stages with that of a known policy with full or partial knowledge of the uncertainty. The regret of an MPC controller for unconstrained problems with process noise has been analyzed in \citep{yu2020power}, for cases where the terminal cost is related to the solution of the Riccati equation. Other related work includes \citep{muthirayan2021online} and  \citep{wabersich2020bayesian}. Although the major challenge addressed by these works comes from uncertainty or partial information, the results echo our conclusion that much of the credit for good performance goes to the MPC structure itself. This provides an alternative perspective for our work and may point to a direction of extensions of the results reported here.

In summary, this paper makes the following contributions:
\begin{itemize}
    \item[(1)] We derive performance bounds for MPC applied to unconstrained and constrained LQR problems, as well as nonlinear systems;
    \item[(2)] We propose a new design for the terminal ingredients, which likely leads to a larger region where MPC is feasible;
    \item[(3)] We conduct numerical studies that verify our theoretical analysis.
\end{itemize}

\section{Performance Bounds for the Linear Quadratic Problem}\label{sec:main}
\subsection{Preliminaries}
We consider discrete-time optimal control problems with system dynamics
\begin{equation}
    \label{eq:dynamics}
    x_{k+1}=Ax_k+Bu_k,\quad k=0,1,\dots,
\end{equation}
and $k$-th stage cost defined
\begin{equation}
    \label{eq:stage_cost}
    x_k'Qx_k+u_k'Ru_k,
\end{equation}
where prime denotes transposition, $R>0$,\footnote{If $X$ and $Y$ are symmetric positive semidefinite matrices, then the notation $X>Y$ ($X\geq Y$) means that the matrix $X-Y$ is positive (semi)definite.} $Q\geq 0$, $(A,B)$ stablizable, and $(A,\sqrt{Q})$ detectable. The optimal feedback law [with respect to cost \eqref{eq:stage_cost} accumulated over an infinite number of stages] is then given by
\begin{equation}
    \label{eq:optimal_control}
    L^*x=-(B'K^*B+R)^{-1}B'K^*Ax,
\end{equation}
where $K^*$ is the unique positive definite solution fulfilling the algebraic equation
\begin{equation}
    \label{eq:riccati}
    K^*=A'\big(K^*-K^*B(B'K^*B+R)^{-1}B'K^*\big)A+Q.
\end{equation}
For this problem, we introduce the Bellman operator $F$ defined for some symmetric matrices as
\begin{equation}
    \label{eq:riccati_op}
    F(K)=A'\big(K-KB(B'KB+R)^{-1}B'K\big)A+Q.
\end{equation}
The $\ell$-fold composition of the operator $F$ is denoted by $F^\ell$. In addition, given a policy $\mu(x)=Lx$, we define the $L$-Bellman operator as
\begin{equation}
    \label{eq:l_riccati_op}
    F_L(K)=(A+BL)'K(A+BL)+Q+L'RL.
\end{equation}
Similarly, the $\ell$-fold composition of the operator $F_L$ is denoted by $F^\ell_L$. A matrix $A$ will be called a \emph{stability matrix} if all of its eigenvalues are less than one in absolute value. Let $\mathbb{S}^n$ denote the space of $n\times n$ real symmetric matrices. For the operators $F$ and $F_L$, the following classical results hold.
\begin{prop}\label{prop:classical_f}
    \begin{itemize}
        \item[(a)] Let $K_1,K_2\in \mathbb{S}^n$ so that $0\leq K_1\leq K_2$. Then we have that $F_L(K_1)\leq F_L(K_2)$ and $F(K_1)\leq F(K_2)$.
        \item[(b)] The equation $K=F(K)$ admits a unique positive definite solution $K^*$. Moreover, $F^k(K)\to K^*$ as $k\to\infty$ for $K\in \mathbb{S}^n$ so that $K\geq 0$.
        \item[(c)] Let $L$ be an $m$ by $n$ matrix such that $(A+BL)$ is a stability matrix. Then the equation $K=F_L(K)$ admits a unique positive definite solution $K_L$. Moreover, $F^k_L(K)\to K_L$ as $k\to\infty$ for $K\in \mathbb{S}^n$ so that $K\geq 0$.
    \end{itemize}
\end{prop}

\subsection{Region of Decreasing}
We will focus on a special set of symmetric matrices, which we call the \emph{region of decreasing}. It is defined as
\begin{equation}
    \label{eq:dos}
    \mathcal{D}=\big\{K\in\mathbb{S}^n\,\big|\,F(K)\leq K \big\}.
\end{equation}
This set of matrices plays an important role in the design of model predictive control, which is a consequence of the definition. The following result can be proved by applying Theorems~10.4 and 10.6 from \citep{bitmead1991riccati}.
\begin{prop}
    \label{prop:rod}
    For all $K\in \mathcal{D}$, $F^k(K)\in \mathcal{D}$ for all $k$; moreover, $K\geq K^*\in \mathcal{D}$.
\end{prop}

Given a matrix $K\in \mathcal{D}$ and a positive integer $\ell$, we are interested in solving the problem via the following scheme:
\begin{equation}\label{eq:unconstrain_mpc}
    \begin{aligned}
        \min_{\{u_k\}_{k=0}^{\ell-1}}& \quad x_{\ell}'Kx_\ell+\sum_{k=0}^{\ell-1} x_k'Qx_k+u_k'Ru_k\\
\mathrm{s.\,t.} & \quad x_{k+1}=Ax_k+Bu_k,\,k = 0,...,\ell-1,\\
	& \quad  x_0 = x
    \end{aligned}
\end{equation}
for all $x$. If the minimum is attained at $(\Tilde{u}_0,\Tilde{u}_1,\dots,\Tilde{u}_{\ell-1})$, the scheme defines a policy $\Tilde{\mu}$ by setting $\Tilde{\mu}(x)=\Tilde{u}_0$; finite horizon LQR intuition shows that there exists a $\Tilde{L}$ such that $\Tilde{\mu}(x)=\Tilde{L}x$. If we apply the Bellman operator and $L$-Bellman operator introduced in \eqref{eq:riccati_op} and \eqref{eq:l_riccati_op} respectively, the policy $\Tilde{\mu}(x)=\Tilde{L}x$ can be defined as
\begin{equation}
\label{eq:sub_l}
    F_{\Tilde{L}}\big(F^{\ell-1}(K)\big)=F^\ell(K).
\end{equation}
Under our assumption, for $K\in \mathcal{D}$, the resulting closed-loop system $(A+B\Tilde{L})$ is stable. This is proved in \citep[Theorem~10.19]{bitmead1991riccati}, and its generalization and connection to performance bound is given in \citep[Section~3.3]{bertsekas2022lessons}. Then there exists a unique positive semidefinite solution to the equation $F_{\Tilde{L}}(K)=K$, which we denote as $K_{\Tilde{L}}$. In the rest of this section, we investigate the suboptimality of the scheme \eqref{eq:unconstrain_mpc}, which is measured through the difference between the infinite horizon cost obtained by the resulting policy $\Tilde{L}$ (encoded through the matrix $K_{\Tilde{L}}$) and the optimal infinite horizon cost (encoded through the matrix $K^*$).

To investigate the properties of $K$, we introduce some additional notation. $\Vert M\Vert$ stands for the induced $2$-norm of the matrix $M$. Following \citep[p.~382]{bertsekas2022abstract} we denote by $\Vert M\Vert_s$ a special weighted Euclidean norm such that $\Vert A+B\Tilde{L}\Vert_s=\sqrt{\rho}<1$ for some $\rho\in(0,1)$. For those two norms, there exist some positive constants $c_1$ and $c_2$ such that
\begin{equation}
    \label{eq:norm_eq}
    c_1\Vert M\Vert\leq \Vert M\Vert_s\leq c_2\Vert M\Vert
\end{equation}
holds for all $n$ by $n$ matrices $M$. An analytical formula for computing $c_1$ and $c_2$ is given in the Appendix~\ref{app:norm_constant}. Both norms have the submultiplicative property. In addition, we define positive integers $\alpha$ and $\beta_\ell$ by $\alpha=\min\{\Vert A+BL^*\Vert^2,1\}$, and $\beta_\ell=\min\{\Vert (A+BL^*)^{\ell-1}\Vert^2,1\}$, where $\ell$ is the horizon length in \eqref{eq:unconstrain_mpc}. Then by the submultiplicative property of the norm, we have that $\beta_\ell\leq \alpha^{\ell-1}.$ The norm $\Vert\cdot\Vert$ also connects to positive semidefiniteness through the following lemma.

\begin{lem}
\label{lma:monotone_norm}
    Let $K_1$ and $K_2$ be symmetric matrices such that $K_1\geq K_2\geq 0$. Then $\Vert K_1\Vert\geq \Vert K_2\Vert$.
\end{lem}

Now we are ready to state our first result.
\begin{prop}\label{prop:sor} 
    For all $K\in \mathcal{D}$ and all integers $i$,
    \begin{equation}
        \label{eq:monotone_1}
        \Vert F^i(K)-K^*\Vert\leq \alpha^i\Vert K-K^*\Vert.
    \end{equation}
    Moreover, if $j$ satisfies $i = (\ell-1)j$, then 
    \begin{equation}
        \label{eq:monotone_2}
        \Vert F^i(K)-K^*\Vert\leq \beta_\ell^j\Vert K-K^*\Vert\leq \alpha^i\Vert K-K^*\Vert.
    \end{equation}
\end{prop}
\begin{pf}
    We first show that
    \begin{equation}
    \label{eq:monotone_bound1}
        \Vert F^i(K)-K^*\Vert\leq \Vert K-K^*\Vert
    \end{equation}
    Since $K\in\mathcal{D}$, we have $F^i(K)\leq K$ and $F^i(K)\geq K^*$ for all $i$. Then 
    $$0\leq F^i(K)-K^*\leq K-K^*.$$
    Taking norm on both sides and applying Lemma~\ref{lma:monotone_norm} gives \eqref{eq:monotone_bound1}. 
    
    By the definition of $F$, it is clear that $F^i(K)\leq F_{L^*}^i(K).$ Combining with the inequality $F^i(K)\geq K^*$ yields
    $$0\leq F^i(K)-K^*\leq F_{L^*}^i(K)-K^*=F_{L^*}^i(K)-F_{L^*}^i(K^*).$$
    Taking norm on both sides and applying Lemma~\ref{lma:monotone_norm}, we have that
    \begin{align*}
        \Vert F^i(K)-K^*\Vert\leq & \Vert F_{L^*}^i(K)-F_{L^*}^i(K^*)\Vert\\
        \leq &\big\Vert\big((A+BL^*)^i\big)' (K-K^*)(A+BL^*)^i\big\Vert\\
        \leq &\Vert(A+BL^*)^i\Vert^2\Vert K-K^*\Vert\\
        \leq & \Vert(A+BL^*)\Vert^{2i}\Vert K-K^*\Vert.
    \end{align*}
    Combing with \eqref{eq:monotone_bound1} yields \eqref{eq:monotone_1}. The inequality \eqref{eq:monotone_2} can be obtained similarly.
\end{pf}

\subsection{Performance Bound via Contraction}\label{sec:contraction}
From the above results we can obtain a performance bound for $K_{\Tilde{L}}$ by exploiting the contraction properties of the operators $F$ and $F_{\Tilde{L}}$.

\begin{prop}\label{prop:con_bound}
    Let $K\in \mathcal{D}$ and $\Tilde{\mu}(x)=\Tilde{L}x$ defined as in \eqref{eq:sub_l}. Then we have that
    \begin{equation}
        \label{eq:lqr_b_contraction}
        \Vert K_{\Tilde{L}} -K^*\Vert\leq  \frac{c_2}{c_1(1-\rho)}(\rho  +\frac{c_2}{c_1}\alpha)\beta_\ell \Vert K-K^*\Vert.
    \end{equation}
\end{prop}

\begin{pf}
Denote as $\overline{K}$ the matrix $F^{\ell-1}(K)$. By Prop.~\ref{prop:sor} with $i=\ell-1$, we have
\begin{equation}
    \label{eq:con_bound_1}
    \Vert \overline{K}-K^*\Vert\leq \beta_\ell \Vert K-K^*\Vert.
\end{equation}
In addition, combining Prop.~\ref{prop:sor} with $i=1$ and the norm equivalence relation \eqref{eq:norm_eq} yields 
\begin{align}
\label{eq:con_bound_2}
    \Vert F(\overline{K})-K^*\Vert_s\leq &c_2\Vert F(\overline{K})-K^*\Vert\nonumber\\
    =&c_2\Vert F(\overline{K})-F(K^*)\Vert\nonumber\\
    \leq &c_2\alpha \Vert \overline{K}-K^*\Vert\leq \frac{c_2}{c_1}\alpha \Vert \overline{K}-K^*\Vert_s,
\end{align}
where the second inequality is due to contraction property of the operator $F$.

We then proceed by showing that 
\begin{equation}
\label{eq:con_bound_3}
    \Vert K_{\Tilde{L}} -K^*\Vert_s\leq \frac{1}{1-\rho}\Vert F_{\Tilde{L}}(K^*)-K^*\Vert_s.
\end{equation}
We use the triangle inequality to write for every $k$,
\begin{align*}
    \Vert F_{\Tilde{L}}^k(K^*)-K^*\Vert_s\leq& \sum_{i=1}^k\Vert F_{\Tilde{L}}^{i}(K^*)-F_{\Tilde{L}}^{i-1}(K^*)\Vert_s\\
    \leq& \sum_{i=1}^k\rho^{i-1}\Vert F_{\Tilde{L}}(K^*)-K^*\Vert_s.
\end{align*}
Taking the limit as $k\to \infty$ and in view of the convergence $F_{\Tilde{L}}^k(K^*)\to  K_{\Tilde{L}}$, we obtain the inequality \eqref{eq:con_bound_3}.

For the value $\Vert F_{\Tilde{L}}(K^*)-K^*\Vert_s$, by using the triangular inequality, and the definition of $\Tilde{L}$ in \eqref{eq:sub_l}, we have
\begin{align}
\label{eq:con_bound_4}
    \Vert F_{\Tilde{L}}(K^*)-K^*\Vert_s\leq& \Vert F_{\Tilde{L}}(K^*)-F_{\Tilde{L}}(\overline{K})\Vert_s+\nonumber\\
    &\Vert F_{\Tilde{L}}(\overline{K})-F(\overline{K})\Vert_s+\Vert F(\overline{K})-K^*\Vert_s\nonumber\\
    =&\Vert F_{\Tilde{L}}(K^*)-F_{\Tilde{L}}(\overline{K})\Vert_s+\Vert F(\overline{K})-K^*\Vert_s\nonumber\\
    \leq &\rho \Vert \overline{K}-K^*\Vert_s +\frac{c_2}{c_1}\alpha \Vert \overline{K}-K^*\Vert_s\nonumber\\
    =&(\rho  +\frac{c_2}{c_1}\alpha)\Vert \overline{K}-K^*\Vert_s,
\end{align}
where in the last step, we use the inequality \eqref{eq:con_bound_2}. Combining the \eqref{eq:con_bound_4} to \eqref{eq:con_bound_3} gives
$$\Vert K_{\Tilde{L}} -K^*\Vert_s\leq \frac{1}{1-\rho}(\rho  +\frac{c_2}{c_1}\alpha)\Vert \overline{K}-K^*\Vert_s.$$
Applying $\Vert K_{\Tilde{L}} -K^*\Vert\leq (1/c_1)\Vert K_{\Tilde{L}} -K^*\Vert_s$ on the left and $\Vert K_{\Tilde{L}} -K^*\Vert_s\leq c_2 \Vert K_{\Tilde{L}} -K^*\Vert$ on the right, we have
$$\Vert K_{\Tilde{L}} -K^*\Vert\leq  \frac{c_2}{(1-\rho)c_1}(\rho  +\frac{c_2}{c_1}\alpha)\Vert \overline{K}-K^*\Vert.$$
In view of the inequality \eqref{eq:con_bound_1}, we get the desired result.
\end{pf}

\begin{rem}
Our proof has followed closely the approach for establishing the classical error bound for the MDP, cf. \citep[Prop.2.2,1]{bertsekas2022abstract}. In particular, if $c_1=c_2$ and $\alpha=\rho$, then the bound established above can be obtained via applying Prop.~2.1.1(e) and Prop.~2.2.1 in \citep[Prop.2.2,1]{bertsekas2022abstract}, and the inequality \eqref{eq:lqr_b_contraction} reduces to the known result in Prop.~2.2.1. This bound indicates that by selecting $K$ close to $K^*$, the MPC policy obtained by solving \eqref{eq:unconstrain_mpc} will perform close to the optimal policy in the infinite horizon.
\end{rem}

\subsection{Performance Bounds via Monotonicity and Newton Step Interpretation}
Alternative bounds can be obtained based on monotonicity and the Newton step interpretation. We first derive the bound that relies solely on the monotonicity, and then provide an additional bound that is obtained by combining the two properties.

\begin{prop}\label{prop:monotone_bound}
    Let $K\in \mathcal{D}$ and $\Tilde{\mu}(x)=\Tilde{L}x$ defined as in \eqref{eq:sub_l}. Then we have that
    \begin{equation}
        \label{eq:lqr_b_monotone}
        \Vert K_{\Tilde{L}} -K^*\Vert\leq  \alpha\beta_\ell\Vert K-K^*\Vert.
    \end{equation}
\end{prop}
\begin{pf}
Recall that the policy $\Tilde{L}$ is defined as
$$F_{\Tilde{L}}\big(F^{\ell-1}(K)\big)=F^\ell(K).$$
Since $K\in \mathcal{D}$, due to the monotonicity property, we have that $F^\ell(K)\leq F^{\ell-1}(K)$. Therefore, applying on both sides $F_{\Tilde{L}}$ and using the definition of $\Tilde{L}$, we have
$$F_{\Tilde{L}}\big(F^{\ell}(K)\big)\leq F_{\Tilde{L}}\big(F^{\ell-1}(K)\big)=F^{\ell}(K).$$
Applying repeatedly on both sides $F_{\Tilde{L}}$ and in view of its monotonicity property, we have $F_{\Tilde{L}}^k\big(F^{\ell}(K)\big)\leq F^{\ell}(K)$ $\forall k>1.$ Taking the limit as $k\to\infty$ and we have $K_{\Tilde{L}}\leq F^{\ell}(K).$ 
On the other hand, we have $K^*\leq K_{\Tilde{L}}$, which implies that
\begin{align*}
    \Vert K_{\Tilde{L}}-K^*\Vert\leq &\Vert F^{\ell}(K)-K^*\Vert=\big\Vert F\big(F^{\ell-1}(K)\big)-K^*\big\Vert\\
    \leq &\alpha  \Vert F^{\ell-1}(K)-K^*\Vert\leq \alpha\beta_\ell  \Vert K-K^*\Vert,
\end{align*}
where the second to last inequality is due to \eqref{eq:monotone_1}, and the last inequality is due to \eqref{eq:monotone_2}.
\end{pf}

\begin{rem}
    Note that by definition $c_2\geq c_1$, and focusing on the last term of \eqref{eq:lqr_b_contraction}, we see that the bound in Prop.~\ref{prop:monotone_bound} is always tighter than the one in Prop.~\ref{prop:con_bound}. This is consistent with what is noted in \citep[Appendix A.3]{bertsekas2022lessons}. Still, the result given in Prop.~\ref{prop:con_bound} resembles the classical performance bound in MDP, and holds true for $K$ where $K\not\in \mathcal{D}$, while the bound in Prop.~\ref{prop:monotone_bound} only holds for $K\in\mathcal{D}$.
\end{rem}

Apart from the monotonicity properties of the operators $F$ and $F_{\Tilde{L}}$, the Newton step interpretation of computing $K_{\Tilde{L}}$ can also be brought to bear for establishing another bound on the suboptimality of $K_{\Tilde{L}}$. For completeness, a slight generalization of the following lemma adopted from \citep[Theorem~2]{hewer1971iterative} and its proof are provided in the Appendix~\ref{app:quadratic}.  
\begin{lem}
\label{lma:newton}
    Let $K\in \mathcal{D}$ and $\Tilde{\mu}(x)=\Tilde{L}x$ defined as in \eqref{eq:sub_l}. Then there exists $\gamma>0$ such that
    \begin{equation}
        \label{eq:lqr_b_newton_tight}
        \Vert K_{\Tilde{L}} -K^*\Vert\leq  \gamma\Vert F^{\ell-1}(K)-K^*\Vert^2.
    \end{equation}
\end{lem}

An explicit formula for $\gamma$ is provided in Appendix~\ref{app:quadratic}. Combining Lemma~\ref{lma:newton} and Prop.~\ref{prop:sor} leads to the following bound.
\begin{prop}\label{prop:newton}
    Let $K\in \mathcal{D}$ and $\Tilde{\mu}(x)=\Tilde{L}x$ defined as in \eqref{eq:sub_l}. Then there exists $\gamma>0$ such that
    \begin{equation}
        \label{eq:lqr_b_newton}
        \Vert K_{\Tilde{L}} -K^*\Vert\leq  \gamma\beta_\ell^2\Vert K-K^*\Vert^2.
    \end{equation}
\end{prop}

\begin{rem}
    The quadratic term $\Vert K-K^*\Vert^2$ in \eqref{eq:lqr_b_newton} manifests the nature of scheme \eqref{eq:unconstrain_mpc} as a step of Newton's method for solving \eqref{eq:riccati}. In addition, if $\gamma\beta_\ell\Vert K-K^*\Vert/\alpha<1$, then the bound based on Newton step interpretation is tighter than the bound provided by Prop.~\ref{prop:monotone_bound}. This is likely to happen if the value $\ell$ is large, which in turn makes the value $\beta_\ell$ small.  
\end{rem}

\begin{rem}
    Since the bound is built upon Lemma~\ref{lma:newton} and Prop.~\ref{prop:sor}, it is more conservative than the bound given in Lemma~\ref{lma:newton}. Still, it connects the matrix $K$ (which is a design parameter), the optimal policy $\mu^*$ (through the parameter $\beta_\ell$) and the performance of $K_{\Tilde{L}}$.     
\end{rem}

\section{Performance Bound for Constrained Problems}\label{sec:mpc}
Let us now consider problems involving both state and control constraints. For the same stationary dynamics \eqref{eq:dynamics} and stage cost \eqref{eq:stage_cost}, there are state constraint $\hat{X}\subset X=\Re^n$ and control constraint $U\subset \Re^m$. We assume that both $\hat{X}$ and $U$ are compact and convex, and contain the origin in their interior.

The problem can be modeled as an optimal control problem involving stationary dynamics \eqref{eq:dynamics} and stage cost 
\begin{equation}
    \label{eq:stage_cost_constraint}
    g(x_k,u_k)=x_k'Qx_k+u_k'Ru_k+\delta_{\hat{X}}(x_k),
\end{equation}
where $\delta_{\hat{X}}$ is an indicator function that maps $x_k$ to $0$ if $x_k\in \hat{X}$, and $\infty$ otherwise. We consider \emph{stationary} policies $\mu$, which are functions mapping $X$ to $U$.

The \emph{cost function} of a policy $\mu$, denoted by $J_\mu$, maps $X$ to $[0,\infty]$, and is defined at any initial state $x_0 \in X$, as
\begin{equation}
\label{eq:pi_cost_function}
    J_\mu(x_0)= \sum_{k=0}^\infty g(x_k,\mu(x_k)),
\end{equation}
subject to $x_{k+1}=Ax_k+B\mu(x_k)$, $k=0,\,1,\,\dots$. The optimal cost function $J^*$ is defined pointwise as
\begin{equation}
    J^*(x_0)=\inf_{\substack{u_k\in U,\ k=0,1,\ldots\\x_{k+1}=Ax_k+Bu_k,\ k=0,1,\ldots}}\sum_{k=0}^\infty  g(x_k,u_k).
\end{equation}
A stationary policy $\mu^*$ is called \emph{optimal} if
\begin{equation*}
    J_{\mu^*}(x)=J^*(x),\quad \forall x\in X.
\end{equation*}
It has the property that
\begin{equation}
\label{eq:mu_star}
    \mu^*(x)\in \min_{u\in U} \{g(x,u)+J^*(Ax+Bu)\},
\end{equation}
if the minimum in \eqref{eq:mu_star} can be attained. For the problem considered here, it may be shown that there is a stationary optimal policy. A brief discussion regarding the existence of a stationary optimal policy is provided in Appendix~\ref{app:existence}.

Similar to the scheme in \eqref{eq:unconstrain_mpc}, we may apply MPC of the following form to solve the problem 
\begin{equation}\label{eq:constrain_mpc}
    \begin{aligned}
        \min_{\{u_k\}_{k=0}^{\ell-1}}& \quad x_{\ell}'Kx_\ell+\sum_{k=0}^{\ell-1} x_k'Qx_k+u_k'Ru_k\\
\mathrm{s.\,t.} & \quad x_{k+1}=Ax_k+Bu_k,\,k = 0,...,\ell-1,\\
& \quad x_k\in \hat{X},\,k = 0,...,\ell-1,\\
& \quad u_k\in U,\,k = 0,...,\ell-1,\\
	& \quad  x_\ell\in S\subset \hat{X},\quad  x_0 = x,
    \end{aligned}
\end{equation}
where $K$ is a suitably designed matrix, and $S$ is a suitably designed set. If the minimum of \eqref{eq:constrain_mpc} is attained at $(\Tilde{u}_0,\Tilde{u}_1,\dots,\Tilde{u}_{\ell-1})$, then the scheme defines a suboptimal policy by setting $\Tilde{\mu}(x)=\Tilde{u}_0$. As in the unconstrained case, we would like to investigate the suboptimality of $\Tilde{\mu}$, depending on the choice of $K$ and $S$, where suboptimality of the policy $\Tilde{\mu}$ is measured by $J_{\Tilde{\mu}}-J^*$.

Due the presence of state constraint $\hat{X}$ and control constraint $U$, the analysis given in Section~\ref{sec:main} does not apply as it is. Still, by suitably modifying the relevant definitions, a similar qualitative analysis remains valid. In particular, we denote by $\mathcal{E}^+(X)$ the set of all functions $J:X\mapsto[0,\infty]$. A mapping that plays the key role in our development is the \emph{Bellman operator} $T:\mathcal{E}^+(X)\mapsto\mathcal{E}^+(X)$, defined pointwise as
\begin{equation}
    \label{eq:bellman_op}
    (TJ)(x)=\inf_{u\in U} \{g(x,u)+J(Ax+Bu)\}.
\end{equation}
This operator is well-posed in view of the nonnegativity \eqref{eq:stage_cost_constraint} of the stage cost. In addition, we denote as $T^\ell$ the $\ell$-fold composition of $T$ with itself, with the convention that $T^0J=J$. For every fixed policy $\mu$, we also introduce the $\mu$-operator $T_\mu:\mathcal{E}^+(X)\mapsto\mathcal{E}^+(X)$, which is defined pointwise as
\begin{equation}
    \label{eq:mu_operator}
    (T_\mu J)(x)=g\big(x,\mu(x)\big)+J\big(Ax+B\mu(x)\big).
\end{equation}
The operators $T$ and $T_\mu$ are generalizations of $F$ and $F_L$. Their relations are extensively discussed in \citep[Chapter~4]{bertsekas2022lessons}.

Given a positive semidifnite matrix $K$ and a set $S$, we can define a function $J\in \mathcal{E}^+(X)$ associated with $K$ and $S$ as 
\begin{equation}
    \label{eq:terminal_j}
    J(x)=x'Kx+\delta_S(x).
\end{equation}
Using the operator notations \eqref{eq:bellman_op}, \eqref{eq:mu_operator}, as well as the function $J$ defined in \eqref{eq:terminal_j}, the policy $\Tilde{\mu}$ defined by the scheme \eqref{eq:constrain_mpc} can be written succinctly as
\begin{equation}
    \label{eq:mu_tilde}
    \big(T_{\Tilde{\mu}}(T^{\ell-1}J)\big)(x)=(T^\ell J)(x), \quad x\in X,
\end{equation}
which is a generalization of \eqref{eq:sub_l}.

For the constrained problem, we can also define the region of decreasing, denoted as $\mathcal{D}(X)$. It is given as a set of functions as
\begin{equation}
    \label{eq:rod}
    \mathcal{D}(X)=\{J\in \mathcal{E}^+(X)\,|\,(TJ)(x)\leq J(x),\,\forall x\in X\}.
\end{equation}

\subsection{Performance Bounds via Monotonicity and Newton Step Interpretation}

For the constrained problem, a result based on contraction, as in Section~\ref{sec:contraction} may not be applicable. However, for suitably designed function $J$, results that parallel Prop.~\ref{prop:monotone_bound} can be established. The proof of the result is deferred to Appendix~\ref{app:pf_mpc_b}.
\begin{prop}
    \label{prop:monotone_bound_con}
    Let $J\in \mathcal{D}(X)$ and $\Tilde{\mu}$ defined via \eqref{eq:mu_tilde}. Then
    \begin{equation}
        \label{eq:mpc_b_1}
        J_{\Tilde{\mu}}(x_0)-J^*(x_0)\leq (T^\ell J)(x_0)-J^*(x_0)\leq J(x_{\ell})-J^*(x_{\ell}),
    \end{equation}
    for all $x_0\in X$ with $x_{k+1}=f\big(x_k,\mu^*(x_k)\big)$, $k=0,\dots,\ell-1$.
    In addition, if $\hat{J}\in \mathcal{E}^+(X)$ such that $\hat{J}(x_0)\leq J^*(x_0)$ for all $x_0\in X$, then
    \begin{equation}
        \label{eq:mpc_b_2}
        J_{\Tilde{\mu}}(x_0)-J^*(x_0)\leq  (T_{\mu^*}^\ell J)(x_0)-(T_{\mu^*}^\ell \hat{J})(x_0)=J(x_{\ell})-\hat{J}(x_{\ell}),
    \end{equation}
    for all $x_0\in X$ with $x_{k+1}=f\big(x_k,\mu^*(x_k)\big)$, $k=0,\dots,\ell-1$.
\end{prop}

\begin{rem}
    A close examination of the proof shows that the above bounds hold well beyond the constrained linear systems discussed here. If a stationary optimal policy $\mu^*$ exists, the performance bounds hold for the case where control constraints are state dependent, i.e., $u\in U(x)$, or problems involving nonlinear, as well as hybrid systems. In fact, the tools applied here have been used to analyze MPC for hybrid systems in \citep{baoti2006constrained}. 
\end{rem}

A Newton's step interpretation of $J_{\Tilde{\mu}}$ similar to that of Prop.~\ref{prop:newton} can also be established. We refer to \citep[Chapter~5]{bertsekas2022lessons} for details.

\subsection{New Designs for the Terminal Costs and Constraints}
For the MPC scheme of the form \eqref{eq:constrain_mpc}, a common choice for $K$ has been $K^*$ given in \eqref{eq:riccati}, i.e., the optimal cost function $x'K^*x$ for the unconstrained LQR problem. The rationale for such a choice is that the corresponding $\Tilde{\mu}$ would be optimal when $x$ is near the origin, thus the performance of $\Tilde{\mu}$ over all feasible $x$ should also be near optimal. However, in the presence of control constraints and if the matrix $L^*$ defined in \eqref{eq:optimal_control} is large, the associated terminal constraint $S$ in \eqref{eq:constrain_mpc} can be small, which means the subset of $X$ from where the scheme \eqref{eq:constrain_mpc} is feasible (or equivalently, the support of $J_{\Tilde{\mu}}$) is small.

On the other hand, our analysis suggests that a choice of $K$ that is near $K^*$ may not be necessary. In fact, one may prefer a much larger $K$ if its associated set $S$ can be larger. A plausible choice that may lead to larger $S$ is the solution of the equation
\begin{equation}
    \label{eq:amplify}
    K=A'\big(K-KB(B'KB+\zeta R)^{-1}B'K\big)A+Q,
\end{equation}
where $\zeta\gg 1$. Such a matrix $K$ belongs to the region of decreasing, and it is likely that the corresponding $S$ is enlarged. This is confirmed in our numerical examples.

We note that a large terminal cost is also advocated in \citep{limon2006stability}, where the terminal constraint is removed. Our study here has deployed different analytical tools and focuses on the classical form with both terminal cost and constraint. Still, it suggested that a terminal cost that corresponds to a smaller control would be preferred.

\section{Numerical Studies}
In this section, we provide numerical examples to demonstrate the differences between the bounds established in the preceding sections and the actual performance. It is observed that the actual performance is often much better than the bounds suggest, further corroborating our main point that it is \emph{not} necessary to design the terminal cost close to optimal at the expense of reduced support of $J_{\Tilde{\mu}}$. 

\subsection{Numerical Study Results for LQR}

\begin{exmp}[Bellman curve for a scalar example]\label{eg:scalar_lqr}
In this example, we consider a scalar system. The purpose is to illustrate the performance bound \eqref{eq:lqr_b_monotone} based on monotonicity and the Newton step interpretation, as shown in Fig.~\ref{fig:scalar}. The system parameters are $A=2,\,B=0.5,\,Q=1,\,R=10$. We set $\ell=1$ so that $K=\overline{K}$. In this case, we have $\Vert K_{\Tilde{L}} -K^*\Vert\approx 3.3$, the bounds \eqref{eq:lqr_b_contraction}, \eqref{eq:lqr_b_monotone} and \eqref{eq:lqr_b_newton} are given as $534.5$, $14.4$ and $43.0$ respectively. The interpretation of Newton's step is evident by noting that $K_{\Tilde{L}}$ is obtained by constructing a tangent line of $F(K)$ at $F(\overline{K})$. 
\begin{figure}[htbp]
    \centering
        \includegraphics[ width=\linewidth]{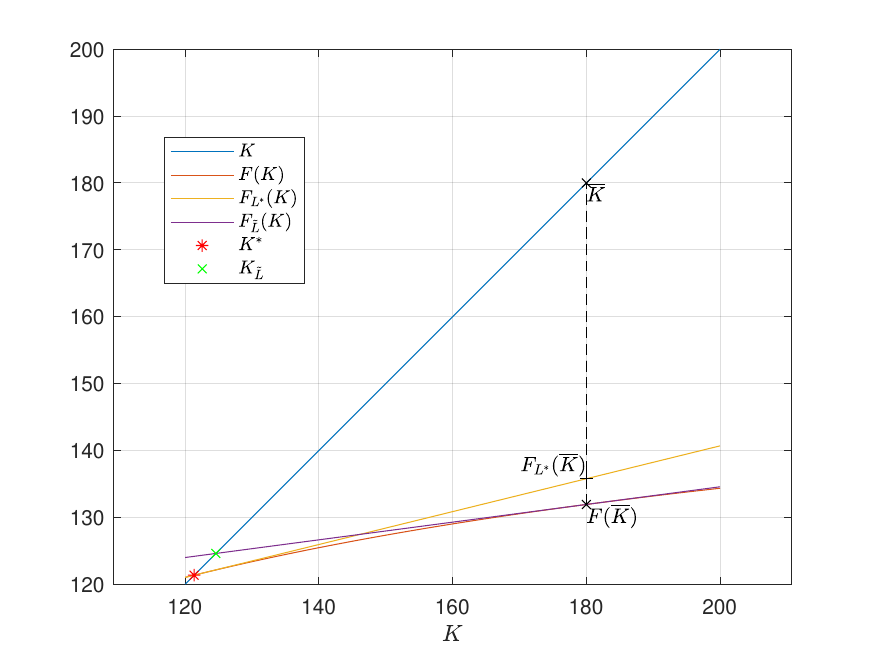}
    \label{fig:scalar}
    \caption{The Bellman function $F$, the $L^*$-Bellman function $F_{L^*}$, the performance bound \eqref{eq:lqr_b_monotone}, as well as the Newton step interpretation of MPC applied to the scalar LQR problem with $\ell=1$. It can be seen that $K^*\leq K_{\Tilde{L}}\leq F(\overline{K})\leq F_{L^*}(\overline{K})$.}
    \end{figure}

\end{exmp}

\begin{exmp}[Conservativeness of bounds]\label{eg:2_4_lqr}
    We consider a two dimensional double integrator and a four dimensional example from \citep{kapasouris1988design}. The problem data of those examples are listed in Appendix~\ref{app:problem_data}. The matrices $K$ in \eqref{eq:unconstrain_mpc} are computed by solving \eqref{eq:amplify} with $\zeta=50$. Their magnitudes as well as their relations to the corresponding $K^*$ are listed in Table~\ref{tb:k}.
\begin{table}[hb]
\begin{center}
\caption{Matrix $K$}\label{tb:k}
\begin{tabular}{ccc}
Problem & $\Vert K\Vert/\Vert K^*\Vert$ & $\Vert K-K^*\Vert$  \\\hline
2-D & $2.5$ & $9.9$\\
4-D & $4.3$ & $486$\\ \hline
\end{tabular}
\end{center}
\end{table}

The actual optimality gap $\Vert K_{\Tilde{L}} -K^*\Vert$ and various performance bounds for those two systems with different $\ell$ are listed in Table~\ref{tb:bounds}. It is clear that the actual performance is much better than what the performance bounds suggest. In addition, the bound \eqref{eq:lqr_b_contraction} is always inferior than \eqref{eq:lqr_b_monotone} and \eqref{eq:lqr_b_newton}, which is consistent with our analysis. 
\begin{table}[hb]
\begin{center}
\caption{Performance bounds with different $\ell$}\label{tb:bounds}
\begin{tabular}{ccccccc}
Problem & $\ell$ & $\Vert K_{\Tilde{L}} -K^*\Vert$ & \eqref{eq:lqr_b_contraction} & \eqref{eq:lqr_b_monotone} & \eqref{eq:lqr_b_newton} \\\hline
2-D & $3$ & $10^{-3}$ & $<10^{10}$ & $9.8$ & $553$\\
2-D & $10$ & $<10^{-13}$ & $<10^{5}$ & $<10^{-4}$ & $<10^{-7}$\\ 
4-D & $3$ &$2.8$ & $10^{52}$ & $486$ & $<10^{10}$\\
4-D & $10$ & $<10^{-3}$ & $>10^{53}$ & $404$ & $<10^{10}$\\ 
4-D & $20$ & $<10^{-7}$ & $<10^{53}$ & $248$ & $<10^{9}$\\
\hline
\end{tabular}
\end{center}
\end{table}

\end{exmp}

\subsection{Numerical Study Results for Constrained LQR}
Here we investigate the performance of the scheme \eqref{eq:constrain_mpc} when both state and control constraints are present. We introduce constraints to the two- and four-dimensional problems studied in Example~\ref{eg:2_4_lqr}. In particular, in Example~\ref{eg:2_mpc}, we demonstrate that with our choice of $K$, the support of $J_{\Tilde{\mu}}$ can be enlarged, while there is hardly any loss of optimality.

\begin{exmp}[Enlarging the feasible region]\label{eg:2_mpc}
We consider the constrained version of the two-dimensional system investigated in Example~\ref{eg:2_4_lqr}. We use the same $K$ as in Example~\ref{eg:2_4_lqr} and set $\ell=3$. The sets $S$ designed according to $K$ and $K^*$ are shown in Fig.~\ref{fig:terminal}. The feasible regions of the MPC with those terminal constraints are illustrated in Fig.~\ref{fig:feasible}. For our choice of $K$, the feasible region practically coincides with the support of $J^*$.

\begin{figure}[htbp]
    \centering
        \includegraphics[ width=\linewidth]{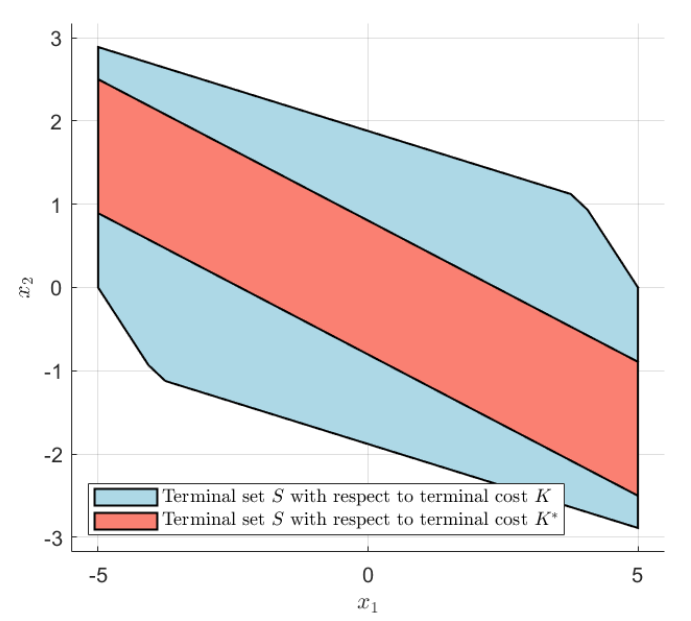}
        \caption{Illustration of terminal sets $S$ designed for $K$ and $K^*$ respectively. It can be seen that for our choice of $K$, which is larger than $K^*$, its corresponding terminal set $S$ is much larger.}
    \label{fig:terminal}
    \end{figure}

    \begin{figure}[htbp]
    \centering
        \includegraphics[ width=\linewidth]{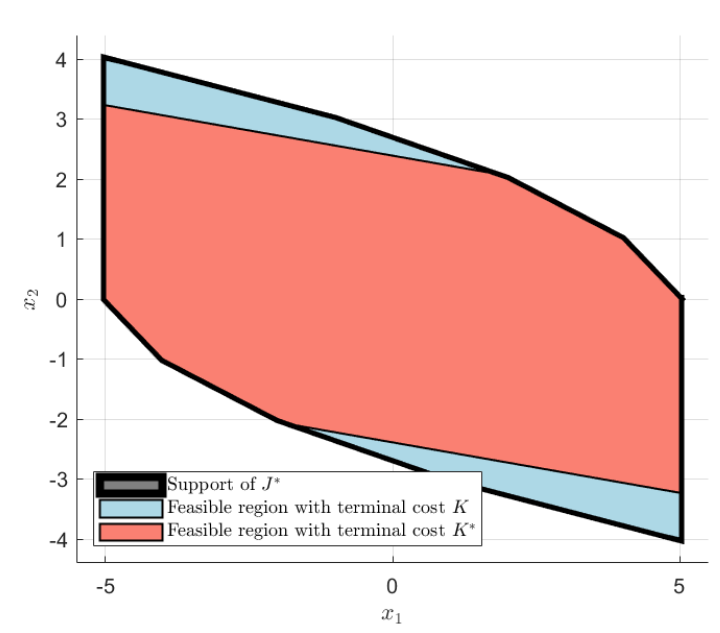}
        \caption{Illustration of feasible region for $\ell=3$ with $K$ and $K^*$ as terminal costs respectively. It can be seen that for our choice of $K$, which is larger than $K^*$, its corresponding feasible set is much larger. In fact, it practically coincides with the support of $J^*$, which is the region from which one can obtain a finite cost.}
    \label{fig:feasible}
    \end{figure}

If we vary the parameter $\zeta$ in \eqref{eq:amplify}, the resulting $K$ would change as well. The following table shows the ratio of the respective terminal set volumes $V_K$ and $V_{K^{\star}}$ for $K$ and $K^{\star}$ with different values of $\zeta$ in \eqref{eq:amplify}.
\begin{table}[hb]
\begin{center}
\caption{Ratio of terminal set volumes with different $\zeta$}\label{tb:volume_ratios}
\begin{tabular}{ccccc}
$\zeta$ & $5$ & $15$ & $25$ & $35$  \\\hline
$V_K/V_{K^{\star}}$ & $1.23$ & $1.56$ &$1.65$ & $1.63$\\
\hline
\end{tabular}
\end{center}
\end{table} 

The suboptimality of the MPC scheme with our choice of $K$ (with $\zeta=50$) and the optimal (approximated in the study by setting $\ell=100$) is given in Fig.~\ref{fig:op_gap}. It can be seen that despite a $K$ that is much larger-than $K^*$, the closed loop system hardly loses any performance.

\begin{figure}[htbp]
    \centering
        \includegraphics[ width=\linewidth]{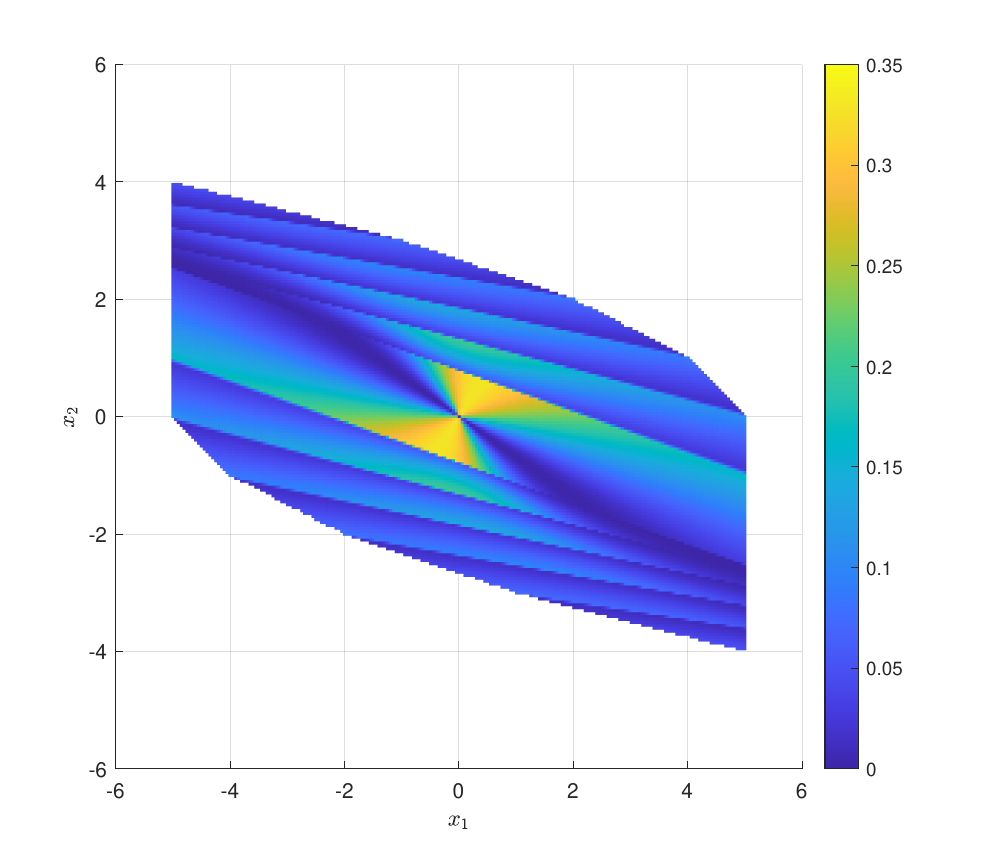}
        \caption{Illustration of the value $|J_{\Tilde{\mu}}(x)-J^*(x)|/|J^*(x)|$ (in percent) with our choice of $K$. It can be seen that there is hardly any loss of optimality. Consistent with our analysis, the largest loss occurs when $x$ is near the origin. Still, a mere $0.35\%$ loss of cost in comparison with optimal is occurred.}
    \label{fig:op_gap}
    \end{figure}

\end{exmp}

\begin{exmp}[Performance loss]\label{eg:4_mpc}
Finally, we investigate the constrained version of the four dimensional system studied in Example~\ref{eg:2_4_lqr}. Using the terminal matrix $K$ applied there and setting $\ell=2$, we show the optimality gap $J_{\Tilde{\mu}}(x)-J^*(x)$ along a trajectory, which is driven by the policy $\Tilde{\mu}$, starting from $x_0$. The result is shown in Fig.~\ref{fig:4d_difference}. To put the values in context, the optimal cost from $x_0$ is $J^*(x_0)\approx 293$. Therefore, the performance is practically optimal.

\begin{figure}[htbp]
    \centering
        \includegraphics[ width=\linewidth]{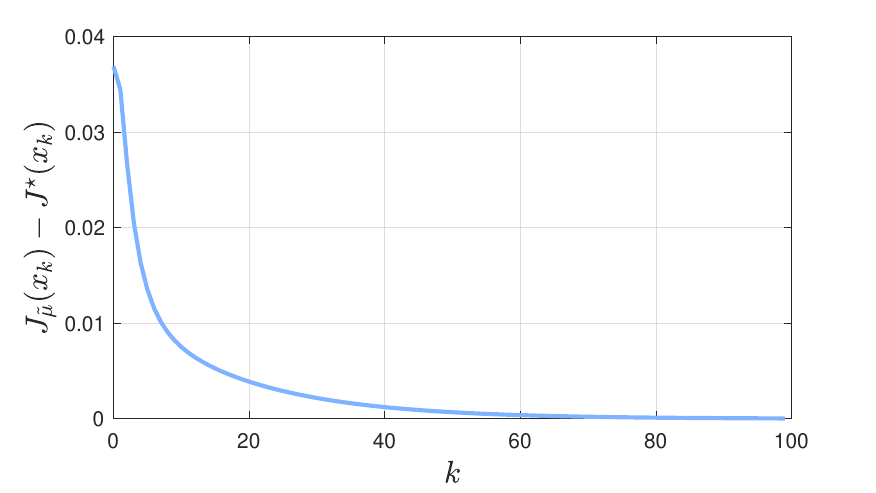}
        \caption{Illustration of the value $J_{\Tilde{\mu}}(x_k)-J^*(x_k)$ where the system is driven under the policy $\Tilde{\mu}$. The optimal cost from $x_0$ is about $293$, which means the policy $\Tilde{\mu}$ is practically optimal.}
    \label{fig:4d_difference}
    \end{figure}
\end{exmp}

\section{Conclusion}
We considered the performance of MPC applied to unconstrained and constrained LQR problems measured by the closed-loop cost accumulated over an infinite horizon. We derived performance bounds that connect the performance of MPC with its terminal cost as well as the true optimal policy of the problems. The derived bounds apply to problems beyond the scope of linear systems and suggest new design of terminal cost and constraint that is not related to the optimal cost of the problem. Numerical studies showed that the new design leads to larger region from which MPC is feasible, while cost little to none performance.

\bibliography{ifacconf}

\appendix
\section{Norm equivalence constants}\label{app:norm_constant}
We provide an analytical formula for the norm equivalence constants $c_1$ and $c_2$ in \eqref{eq:norm_eq}. The results are summarized in the following proposition. A related result is given as \citep[Theorem~5.6.7]{horn2012matrix}.

\begin{prop}
    Let $W$ be an $n$ by $n$ nonsingular matrix and the norm $\Vert M\Vert_s$ for a matrix $M$ is defined as
    $$\Vert M\Vert_s=\max_{x\neq 0}\frac{\Vert Mx\Vert_s}{\Vert x\Vert_s},$$
    where $\Vert x\Vert_s=\Vert Wx\Vert$. Then the norm equivalence condition \eqref{eq:norm_eq} holds with
    $$c_1=\frac{\lambda_n}{\lambda_1},\quad c_2=\frac{\lambda_1}{\lambda_n},$$
    where $\lambda_1$ and $\lambda_n$ are largest and smallest eigenvalues of $W'W$, respectively.
\end{prop}

\begin{pf}
    For $c_2$ note that
    \begin{align*}
        \Vert M\Vert_s&=\max_{x\neq 0}\frac{\Vert Mx\Vert_s}{\Vert x\Vert_s}=\max_{x\neq 0}\frac{\Vert WMx\Vert}{\Vert Wx\Vert}\\
        \leq & \Vert W\Vert\Vert M\Vert\max_{x\neq 0}\frac{\Vert x\Vert}{\Vert Wx\Vert}.
    \end{align*}
    Also note that $\Vert W\Vert=\lambda_1$ and 
    \begin{equation}
    \label{eq:lambda_n}
        \max_{x\neq 0}\frac{\Vert x\Vert}{\Vert Wx\Vert}=\max_{x\neq 0}\sqrt{\frac{x'x}{x'W'Wx}}=\lambda_n,
    \end{equation}
    which gives a formula for $c_2$.

    Similarly, for $c_1$, we have that
    \begin{align*}
        \Vert M\Vert&=\max_{x\neq 0}\frac{\Vert Mx\Vert}{\Vert x\Vert}=\max_{x\neq 0}\frac{\Vert W^{-1}Mx\Vert_s}{\Vert W^{-1}x\Vert_s}\\
        \leq & \Vert W^{-1}\Vert_s\Vert M\Vert_s\max_{x\neq 0}\frac{\Vert x\Vert_s}{\Vert W^{-1}x\Vert_s}.
    \end{align*}
    Moreover
    \begin{equation*}
        \Vert W^{-1}\Vert_s=\max_{x\neq 0}\frac{\Vert WW^{-1}x\Vert}{\Vert Wx\Vert}=\lambda_n,
    \end{equation*}
    as given in \eqref{eq:lambda_n}. In addition, we also have that
    \begin{align*}
        \max_{x\neq 0}\frac{\Vert x\Vert_s}{\Vert W^{-1}x\Vert_s}=&\max_{x\neq 0}\frac{\Vert Wx\Vert}{\Vert WW^{-1}x\Vert}\\
        =&\max_{x\neq 0}\sqrt{\frac{x'W'Wx}{x'x}}=\lambda_1,
    \end{align*}
    which completes the proof.
\end{pf}

\section{Proof for quadratic convergence}\label{app:quadratic}
We first list two useful lemmas. The first one can be found, among others, in \citep[Corollary~2]{lopez2021efficient}, and the second is derived from \citep[Proof of Theorem 1 part~2)]{hewer1971iterative}.
\begin{lem}
    Let $D$ be a stability matrix, then $\sum_{i=0}^\infty\Vert D^i\Vert^2<\infty$.
\end{lem}

\begin{lem}\label{lma:optimality_gap}
Let $\Tilde{L}$ be an $m$ by $n$ matrix such that $\Tilde{D}=A+B\Tilde{L}$ is a stability matrix. Denote as $K_{\Tilde{L}}$ the matrix satisfying the equation $K=F_{\Tilde{L}}(K)$. Let $K^*$ and $L^*$ be the matrices defined by \eqref{eq:riccati} and \eqref{eq:optimal_control}. Then we have that
$$K_{\Tilde{L}}-K^*=\Tilde{D}'(K_{\Tilde{L}}-K^*)\Tilde{D}+(\Tilde{L}-L^*)'(B'K^*B+R)(\Tilde{L}-L^*).$$
\end{lem}

Based on the above equality, we are ready to prove the quadratic convergence result.
\begin{prop}\label{prop:rollout_1_superli}
Given $\overline{K}\in \mathcal{D}$, there exists some $\gamma>0$, such that
$$\Vert K_{\Tilde{L}} -K^*\Vert\leq  \gamma\Vert \overline{K}-K^*\Vert^2,$$
where $\Tilde{L}$ is defined by $F_{\Tilde{L}}(\overline{K})=F(\overline{K})$, and $K_{\Tilde{L}}$ satisfies the equation $K_{\Tilde{L}}=F_{\Tilde{L}}(K_{\Tilde{L}})$.
\end{prop}

\begin{pf}
Due to Lemma~\ref{lma:optimality_gap}, we have that
$$K_{\Tilde{L}}-K^*=\Tilde{D}'( K_{\Tilde{L}}-K)\Tilde{D}+(\Tilde{L}-L^*)'(B'K^*B+R)(\Tilde{L}-L^*).$$
The difference $K_{\Tilde{L}}-K^*$ can be computed as
$$K_{\Tilde{L}}-K^*=\sum_{i=0}^\infty (\Tilde{D}')^i(\Tilde{L}-L^*)'(B'K^*B+R)(\Tilde{L}-L^*)\Tilde{D}^i.$$
Using induced $2$-norm for the matrix, we have
\begin{align}
\label{eq:ineq1}
    &\Vert K_{\Tilde{L}}-K^*\Vert\nonumber\\
    \leq &\sum_{i=0}^\infty\Vert (\Tilde{D}')^i(\Tilde{L}-L^*)'(B'K^*B+R)(\Tilde{L}-L^*)\Tilde{D}^i\Vert\nonumber\\
    \leq &\Vert(\Tilde{L}-L^*)'(B'K^*B+R)(\Tilde{L}-L^*)\Vert\sum_{i=0}^\infty\Vert \Tilde{D}^i\Vert^2\nonumber\\
    \leq &\Vert \Tilde{L}-L^*\Vert^2\Vert B'K^*B+R\Vert\sum_{i=0}^\infty\Vert \Tilde{D}^i\Vert^2,
\end{align}
where $\sum_{i=0}^\infty\Vert \Tilde{D}^i\Vert^2<\infty$ as $\overline{K}\in \mathcal{D}$. 

Now consider the term $\Tilde{L}-L$. By straightforward computation, we have that
\begin{align*}
    \Tilde{L}-L^*=&(B'K^*B+R)^{-1}B'\big((K^*-\overline{K})A+\\
    &(\overline{K}-K^*)B(B'\overline{K}B+R)^{-1}B'\overline{K}A\big).
\end{align*}
Thus, we have that
\begin{equation}
    \label{eq:ineq2}
    \Vert \Tilde{L}-L^*\Vert\leq\eta\Vert K^*-\overline{K}\Vert,
\end{equation}
where $\eta$ is given as
\begin{align*}
    \eta=&\Vert (B'K^*B+R)^{-1}\Vert\cdot\\
    &\big(\Vert B'\Vert \Vert A \Vert+\Vert B'\Vert \Vert B \Vert \Vert (B'\overline{K}B+R)^{-1}\Vert \Vert B'\overline{K}A\Vert\big)
\end{align*}

Combining \eqref{eq:ineq1} and \eqref{eq:ineq2} yields that
\begin{equation}
    \label{eq:ineq3}
    \Vert K_{\Tilde{L}} -K^*\Vert\leq\gamma \Vert \overline{K}-K^*\Vert^2,
\end{equation}
where $\gamma$ is 
$$\gamma = \eta^2\Vert B'K^*B+R\Vert\sum_{i=0}^\infty\Vert \Tilde{D}^i\Vert^2$$

\end{pf}

\section{Existence of a stationary optimal policy}\label{app:existence}
We give a brief discussion on the existence of a stationary optimal policy for the constrained problems considered in Section~\ref{sec:mpc}. We provide the outline according to which one can establish the existence of a stationary optimal policy. Since such a result is peripheral to the main point of this work, the detailed discussion is omitted.

Let $J_0\in \mathcal{E}^+(X)$ be the function that is identically zero. Define the sequence of functions $\{J_k\}$ by
$$J_{k+1}=TJ_k,\quad k=0,1,\dots.$$
For every fixed $x$ and $\lambda\geq0$, consider the sets
$$U_k(x,\lambda)=\{u\in U\,|\,g(x,u)+J_k(Ax+Bu)\leq \lambda\}$$
defined for $k=0,1,\dots$. If there exists some $\Bar{k}$ such that $U_k(x,\lambda)$ is compact for all $k\geq \Bar{k}$, then we can assert the existence of the stationary optimal policy. This condition is given in the form above in \citep{bertsekas1975monotone}. A related condition is derived independently in \citep{schal1975conditions}.

\section{Proof of Prop.~\ref{prop:monotone_bound_con}}\label{app:pf_mpc_b}

We first introduce a useful result, which holds well-beyond the scope considered here \citep[Prop.~5.2]{bertsekas1978stochastic}. It underlies the importance of the region of decreasing introduced in this work, and is closely related to the control Lyapunov function used in typical MPC analysis, see, e.g., \citep[Definition B.39]{rawlings2017model}. The result stated below is adopted from \citep[Prop.~4.3.4]{bertsekas2022abstract}.

\begin{lem}\label{lma:cost_bound}
    Let $J\in \mathcal{E}^+(X)$ and $\mu$ be some stationary policy. Then we have that
    \begin{equation}
        \label{eq:bellman_mu}
        J_\mu(x)=(T_\mu J_\mu)(x),\quad \forall x\in X.
    \end{equation}
    Moreover, if the inequality $(T_\mu J)(x)\leq J(x)$ holds for all $x\in X$, then $J_\mu(x)\leq J(x)$, $\forall x\in X$.
\end{lem}

Now we are ready to prove Prop.~\ref{prop:monotone_bound_con}.
\begin{pf}
    Since $J\in \mathcal{D}(X)$, then $T^{\ell-1}J\in\mathcal{D}(X)$. Therefore, we have 
    $$(T^{\ell}J)(x_0)\leq (T^{\ell-1}J)(x_0),\quad \forall x_0\in X.$$
    Since the policy $\Tilde{\mu}$ is defined via \eqref{eq:mu_tilde}, applying $T_{\Tilde{\mu}}$ on both sides of the above inequality, we have that
    $$\big(T_{\Tilde{\mu}}(T^{\ell}J)\big)(x_0)\leq\big(T_{\Tilde{\mu}}(T^{\ell-1}J)\big)(x_0)=(T^\ell J)(x_0),$$
    where the equality is due to \eqref{eq:mu_tilde}. By Lemma~\ref{lma:cost_bound}, we have $T^\ell J$ as an upper bound of $\Tilde{\mu}$, i.e.,
    $$J_{\Tilde{\mu}}(x_0)\leq (T^\ell J)(x_0),\quad \forall x_0\in X.$$
    Subtracting $J^*(x_0)$ on both sides of the inequality leads to the left inequality of \eqref{eq:mpc_b_1}. 

    The right inequality in \eqref{eq:mpc_b_1} follows from the definition of $T$. In particular, for every $J\in\mathcal{E}^+(X)$, every stationary policy $\mu$ and every integer $\ell$, we have that 
    $$(T^\ell J)(x_0)\leq (T^\ell_\mu J)(x_0),\quad \forall x_0\in X.$$
    The inequality follows by considering the optimal policy $\mu^*$. In particular, we have that
    \begin{align}
        &(T^\ell J)(x_0)-J^*(x_0)\leq (T_{\mu^*}^\ell J)(x_0)-J^*(x_0)\nonumber\\
        =&J(x_\ell)+\sum_{k=0}^{\ell-1}g\big(x_k,\mu^*(x_k)\big)-J^*(x_0),\label{eq:prop13_b1}
    \end{align}
    where $x_{k+1}=f\big(x_k,\mu^*(x_k)\big)$, $k=0,\dots,\ell-1$. In view of \eqref{eq:bellman_mu} and the fact that $J_{\mu^*}(x)=J^*(x)$ for all $x$, we have that
    \begin{equation}
        \label{eq:prop13_b2}
        J^*(x_0)=J^*(x_\ell)+\sum_{k=0}^{\ell-1}g\big(x_k,\mu^*(x_k)\big),
    \end{equation}
    where $x_{k+1}=f\big(x_k,\mu^*(x_k)\big)$, $k=0,\dots,\ell-1$. Combining \eqref{eq:prop13_b1} and \eqref{eq:prop13_b2} together gives the desired result.
    
    The inequalities in \eqref{eq:mpc_b_2} are obtained by noting that $\hat{J}(x)\leq J^*(x)$.
\end{pf}

\section{Problem data for numerical examples}\label{app:problem_data}
For the two dimensional example in Example \ref{eg:2_4_lqr} and in Example \ref{eg:2_mpc} the system matrices are
\begin{equation*}
    A = \begin{bmatrix}
1 & 1 \\
0 & 1
\end{bmatrix}, \quad 
B = \begin{bmatrix}
0\\
1
\end{bmatrix},
\end{equation*}
and the cost matrices are taken to be $R=1$ and the identity matrix of dimension two for $Q$. In Example \ref{eg:2_mpc} the state and input constraints are $|x_1|\leq 5,|x_2|\leq 5$ and $|u|\leq1$. For this example, the matrix $K$ is taken to be the solution of \eqref{eq:amplify} with $\zeta=50$.

For the four dimensional example in Examples \ref{eg:2_4_lqr} and \ref{eg:4_mpc}, the system matrices are taken to be the approximate time-invariant model studied in \citep{kapasouris1988design}. The continuous dynamics are discretized with a zero-order-hold and a sampling time of $50\,\text{ms}$ resulting in the following discrete time matrices
\begin{align*}
    A &= \begin{bmatrix}
0.9993 & -3.0083 & -0.1131 & -1.6081 \\
0 & 0.9862 & 0.0478 & 0 \\
0 & 2.0833 & 1.0089 & 0 \\
0 & 0.0526 & 0.0498 & 1 \\
\end{bmatrix}, \quad \\
B &= \begin{bmatrix}
-0.0804 & -0.6347 \\
-0.0291 & -0.0143\\
-0.8679 &  -0.0917\\
-0.0216 & -0.0022
\end{bmatrix}.
\end{align*}
In Example \ref{eg:2_mpc} the state and input constraints are $|x_2|\leq 0.5$ and $|u_1|\leq25,|u_2|\leq25$. The matrix $K$ for this example has been taken to be the solution of \eqref{eq:amplify} with $\zeta=50$. 

The numerical Examples~\ref{eg:2_mpc} and \ref{eg:4_mpc} are solved in MATLAB with toolbox MPT \citep{kvasnica2004multi} and solver \texttt{quadprog}.

\end{document}